%% file: main.tex
\documentclass{article}
\usepackage{graphicx,url}
\usepackage{epsfig,epsf,listings}
\usepackage{amsmath,amssymb,amsfonts,bbm,mathrsfs}
\usepackage{caption}
\usepackage{float}
\usepackage[version=3]{mhchem}
\usepackage{ulem}
\usepackage{subfigure}
\usepackage{fancyhdr}
\usepackage{epstopdf}
\usepackage[utf8]{inputenc}
\usepackage[utf8]{inputenc}
\usepackage{enumerate}
\usepackage{mathbbol}
\usepackage{dsfont}
\usepackage{multirow,verbatim,color}
\usepackage{algorithm}
\usepackage{algpseudocode}

\newcommand{\K}{\mathcal{K}}

\newcommand{\V}{\mathbf{V}}

\newcommand{\Hb}{\mathbf{H}}
\newcommand{\A}{\mathbf{A}}

\renewcommand{\xi}{\overline{x}}

\renewcommand{\thispagestyle}[2]{}

\fancypagestyle{plain}{
        \fancyhead{}
        \fancyhead[C]{first page center header}
        \fancyfoot{}
        \fancyfoot[C]{first page center footer}
}
\begin{document}
\include{logo}
\begin{abstract}
This study considers using Metropolis-Hastings algorithm for
stochastic simulation of chemical reactions. The proposed method uses SSA (Stochastic Simulation Algorithm) distribution which is a standard method for solving well-stirred chemically reacting systems as a desired distribution. A new numerical solvers based on exponential form of exact and approximate solutions of CME (Chemical Master Equation) is employed for obtaining target and proposal distributions in Metropolis-Hastings algorithm to accelerate the accuracy of the tau-leap method. Samples generated by this technique have the same distribution as SSA and the histogram of samples show it's convergence to SSA.

\paragraph{key words} Metropolis-Hastings, SSA, CME, tau-leap.
\end{abstract}

\section{Introduction}
\label{sect:intro}
In biological systems chemical reactions are modeled stochastically. The system's state (the number of of molecules of each individual species) is described by probability densities describing the quantity of molecules of different species at a given time. The evolution of probabilities through time is described by the chemical master equation (CME) \cite{Gillespie_1977}.

The Stochastic Simulation Algorithm (SSA) first introduced by Gillespie \cite{Gillespie_1977}, is a Monte Carlo approach to sample from the CME. The accuracy of different approaches in simulating stochastic chemical reactions is compared to histogram of samples obtained by SSA. However, SSA has a number of drawbacks such as it simulates one reaction at a time  and therefore it is inefficient for most realistic problems. Alternative approaches have been developed trying to enhance the efficiency of SSA but most of them suffer from accuracy issues. 
The explicit tau-leaping method \cite{Gillespie_2001} iis able to simulate multiple chemical reactions in a pre-selected time step of length $\tau$ by using Poisson random variables \cite{Kurz_1972_SSA}. However, explicit tau-leaping method is numerically unstable for stiff systems  \cite{Cao_2004_stability}. Different implicit tau-leap approaches have been proposed to alleviate the stability issue \cite{Rathinam_2003,Gillespie_2003,Gillespie_2001,Sandu_2013_SSA}. 
Sandu  \cite{Sandu_2013_CME} considers an exact exponential solution to the CME, leading to a solution vector that coincides with the probability of SSA. Several approximation methods to the exact exponential solution as well as approximation to the explicit tau-leap are given in \cite{Azam_Sandu_ApproximateCME}. 

The availability of exact and  approximate probability solutions motivates the use of Markov chain metropolis algorithm to enhance the accuracy of explicit tau-leap method when using large time steps. The proposed method relies on explicit tau-leaping to generate candidate samples. The proposed probability density corresponds to that of tau-leaping \cite{Sandu_2013_CME}, and the target probability density is provided by the CME.   During the Markov process the candidate samples are evaluated based on approximations of target and proposal probability and are either accepted or rejected. The proposed technique requires the computation of a matrix exponential during the Markov process. The dimension of matrix grows with increasing number of species in a reaction system. In order to manage the computational expense of matrix exponentiation efficient approaches based on Krylov \cite{expokit} and rational approximations \cite{pade, cram} are employed. Further computational savings are obtained by exponentiating only a sub-matrix that encapsulates the essential information about the transition of the system from the current to the proposed state.

The paper is organized as follows. Section \ref{sect:mcmc} reviews Monte Carlo approaches, 
and Section \ref{sect:MH_stochastic} discusses the application of Metropolis Hastings algorithm to sample from the probability distribution generated by CME. Computationally efficient methods to accelerate exponentiating the matrix are discussed in Section \ref{sect:mexp}. Numerical experiments carried out in Section \ref{sect:experiment} illustrate the accuracy of the proposed schemes. Conclusions are drawn in Section \ref{sect:conc}.

\section{Markov Chain Monte Carlo methods}
\label{sect:mcmc}
Markov Chain Monte Carlo (MCMC) methods are a class of algorithms to generate samples from desired probability distributions.  A Markov chain is a discrete time stochastic process, i.e., a sequence of random variables (states) ${x_0,x_1, \cdots}$ where the probability of the next state of system depends only on the current state of the system and not on previous ones \cite{Metro_intro}.


\subsection{Metropolis methods}
\label{sect:metropolis}
The Metropolis method is an MCMC process to obtain a sequence of random samples from a desired probability distribution $\pi(x)$, $x \subset X  \subset \mathbb{R}^n$, which is usually complex. 
A Markov chain with state space $X$ and equilibrium distribution $\pi(x)$ is constructed and long runs of the chain are performed \cite{Metropolis_1994}.
The original MCMC algorithm was given by Metropolis et al. \cite{Metropolis_main} and was later modified by Hastings \cite{Metropolis_hasting_main}, with a focus on statistical problems.

A random walk is performed around the current state of the system $x_{t-1}$. A proposal distribution 
$g \left(x^{*}|x_{t-1}\right)$ is used to suggest a candidate $x^{*}$ for the next sample  given the previous sample value $x_{t-1}$. The proposal distribution should be symmetric $ g\left(x_{t-1}|x^{*}\right)=g\left(x^{*}|x_{t-1}\right)$. 
The algorithm works best if the proposal density matches the shape of the target distribution, i.e. $g \left(x_{t-1}|x^{*}\right) \approx \pi(x)$. Proposals $x^{*}$ are accepted or rejected in a manner that leads system toward the region of higher target probability $\pi(x)$ \cite{Metro_history}. Specifically, one computes the target density ratio 
\begin{equation} \label{eq-mc-ratio}
\alpha=\frac{\pi\left(x^{*}\right)}{\pi\left(x_{t-1}\right)}
\end{equation}
and draws a random variable $\zeta \sim uniform(0,1)$. The proposal is accepted or rejected as follows:
\[
x_{t} := \left\{
\begin{array}{lll}
x^{*} & if \quad \zeta < \min \left(1,\alpha\right) & \textnormal{(proposal accepted)},\\
x_{t-1} & otherwise & \textnormal{(proposal rejected)}.
\end{array}
\right.
\]

\subsection{Metropolis Hastings algorithm}
\label{sect:hastings}
The Metropolis Hastings algorithm allows more freedom in the choice of the proposal distribution by relaxing the symmetry constraint \cite{Metro_understand}. 
The acceptance ratio \eqref{eq-mc-ratio} is changed to
\begin{equation} \label{eq-ratio}
\alpha=\alpha_1 \cdot \alpha_2.
\end{equation}
Here $\alpha_{1}$ is the  ratio between the target probabilities of the proposal sample $x^{*}$ and of the previous sample $x_{t-1}$. It can be evaluated by a function $f$ which is an approximation of $\pi$
\begin{equation} \label{eq-ratio1}
\alpha_{1}=\frac{\pi\left(x^{*}\right)}{\pi\left(x_{t-1}\right)
} \approx  \frac{f \left(x^{*}\right)}{f\left(x_{t-1}\right)}
\end{equation}
The ratio $\alpha_{2}$ of the proposal densities of $x^{*}$ conditioned by $x_{t-1}$, and of $x_{t-1}$ conditioned by $x^{*}$ is equal to one if the proposal distribution is symmetric
\begin{equation} \label{eq-ratio2}
\alpha_{2}=\frac{g \left(x_{t-1}|x^{*}\right)}{g \left(x^{*}|x_{t-1}\right)}.
\end{equation}

Convergence of the Markov chain is guaranteed if the properties of detailed balance and ergodicity conditions are fulfilled \cite{Metro_converg}. 
Detailed balance requires that the probability of moving from $x_{t-1}$ is the same as moving from $x^*$.
\[
\pi \left(x_{t-1}\right)g\left(x_{t-1}|x^*\right)=\pi \left(x^*\right)g\left(x^*|x_{t-1}\right)
\]
Ergodicity requires that a chain starting from any state $x_1$ 
will return to $x_1$  if it runs long enough.
In practice, it is not possible to establish with full certainty that a chain has converged \cite{Metro_converg}.

\section{Metropolis Hastings for stochastic simulation}
\label{sect:MH_stochastic}
Here we discuss the application of the Metropolis Hastings algorithm to generate samples from the CME distribution. SSA is currently the standard model for solving well-stirred chemically reacting systems; however, SSA does one reaction at a time that making it slow for real problems. On the other hand, alternative techniques such as explicit and implicit tau-leap methods are faster than SSA but suffer from low accuracy at larger time steps.

In the proposed approach, explicit tau-leap is employed to generate candidate samples. The samples are evaluated based on the acceptance ratio of Metropolis Hastings algorithm. At the end of algorithm, the samples generated by this technique have the same distribution as given by CME, and the histogram of samples converges to the histogram of SSA solutions. 


\subsection{Target distribution}
\label{sect:target}
The target (exact) distribution $\mathcal{P}\left(x,t\right)$ of the state of the chemical system is given by the solution of the CME \cite{Gillespie_1977} 
\begin{equation} \label{cme}
\frac {\partial\mathcal{P}\left(x,t\right)}{\partial t}=\sum_{r=1}^M a_{r}\left(x-v_{r}\right)
\mathcal{P}\left(x-v_{r},t\right)-\\ 
a_0\left(x\right)\mathcal{P}\left(x,t\right)\,.
\end{equation}
Let $Q^i$ is the total possible number of molecules of species $S^i$, $i=1,\dots,N$. The
total number of all possible states of the system is 
\begin{equation}  \label{eq_Q}
Q=\prod_{i=1}^{N}\left(Q^i+1\right).
\end{equation}
CME is a linear ODE on the discrete state space of states $\mathbb{R}^Q$
\begin{equation}  \label{eq_cme_mat}
\mathcal{P}' = A \cdot \mathcal{P}\,, \quad  \mathcal{P}(\bar{t}) = \delta_{\mathcal{I}(\bar{x})}\,, \quad
t \ge \bar{t}\,.
\end{equation}
and has an exact solution: 
\begin{equation}  \label{eq_exact}
\mathcal{P}\left(\bar{t}+\tau\right)=\exp\left(T\, A\right)\cdot \mathcal{P}\left(\bar{t}\right)
= \exp\left(\tau\, \sum_{r=0}^M A_r\right)\cdot  \mathcal{P}\left(\bar{t}\right)\,.
\end{equation}
As explained  in \cite{Sandu_2013_CME,Azam_Sandu_ApproximateCME}, the diagonal matrix  $A_{0}\in \mathbb{R}^{Q \times Q} $ and the Toeplitz matrices $A_{1},\cdots,A_{M}\in \mathbb{R}^{Q \times Q} $ are:
\begin{equation}  \label{eq_sumofexponentexact}
({A_{0}})_{i,j}=\left\{ 
\begin{array}{rl}
-a_{0}\left(x_j\right) & \mbox{if $i=j$} ,\\ 
0 & \mbox{if $i \not =j$} ,
\end{array}
\right. \,, \quad
({A_{r}})_{i,j}=\left\{ 
\begin{array}{rl}
a_{r}(x_j) & \mbox{if $i-j=d_{r}$}, \\ 
0 & \mbox{if $i-j \not =d_{r}$},
\end{array}
\right.
\end{equation}
and their sum $A \in \mathbb{R}^{Q \times Q}$ is
\begin{equation}  \label{eq_exact1}
A = A_{0} +  \dots + A_{M}\,, \quad
A_{i,j}=\left\{ 
\begin{array}{rl}
-a_{0}(x_j) & \mbox{if }i=j\,, \\ 
a_{r}(x_j) & \mbox{if }i-j=d_{r},~ r=1,\cdots,M\,, \\ 
0 & \mbox{otherwise} \,.
\end{array}
\right.
\end{equation}
Here $x_j$ denotes the unique state with state space index $j=\mathcal{I}(x_j)$,
where $\mathcal{I}(x)$ is the state-space index of state $x=[X^1,\, \ldots\,,X^N]$:
\begin{equation}
\begin{array}{rcl} \label{eq_inx}
\mathcal{I}(x) &=& \left(Q^{N-1}+1\right)\cdots \left(Q^1+1\right)\cdot X^N+\cdots  \\ 
&& +\left(Q^2+1\right)\left(Q^1+1\right)\cdot X^3+\left(Q^1+1\right)\cdot X^2+X^1+1. 
\end{array}%
\end{equation}
One firing of reaction $R_{r}$ changes the state from $x$ to $\bar {x}=x-v_{r}$.
The corresponding change in state space index is: 
\[  \label{eq_d}
\begin{array}{l}
\mathcal{I}(x)-\mathcal{I}\left(x-v_{r}\right)=d_{r},   \\ 
d_{r}=\left(Q^{N-1}+1\right)\cdots\left(Q^1+1\right).v_{r}^N+ \cdots \\
+\left(Q^2+1\right)\left(Q^1+1\right).v_{r}^3+\left(Q^1+1\right).v_{r}^2+v_{r}^1.
\end{array}
\]
At the current time $\bar{t}$ the system is in the known state $x(\bar(t))=\bar{x}$ and consequently
the current distribution $\mathcal{P}\left(\bar{t}\right) = \delta_{\mathcal{I}(\bar{x})}$ is equal to one at $\mathcal{I}(\bar{x})$ and is zero everywhere else.
The target distribution in our method is the exact solution \eqref{eq_exact}
\begin{equation}  \label{eq_target}
\pi = \exp\left(\tau\, \sum_{r=0}^M A_r\right)\cdot \delta_{\mathcal{I}(\bar{x})}\,.
\end{equation}
%
\subsection{Proposal distribution}
\label{sect:proposal}
In our algorithm the explicit tau-leap method is employed to generate 
the candidate samples. Sandu \cite{Sandu_2013_CME} shows that the probability distribution generated by the tau-leap method is the solution of a linear approximation of the CME
\begin{equation}  \label{eq_sumofexponent3}
\mathcal{P}\left(\bar{t}+\tau\right)=\exp\left(\tau\, \bar{A}\right)\cdot \mathcal{P}\left(\bar{t}\right)
= \exp\left(\tau\, \sum_{r=0}^M \bar{A}_r\right)\cdot  \mathcal{P}\left(\bar{t}\right)
\end{equation}
where the diagonal matrix  $\bar{A_{0}}\in \mathbb{R}^{Q \times Q} $ and the Toeplitz matrices 
$\bar{A_{1}},...,\bar{A_{M}}\in \mathbb{R}^{Q \times Q} $ are:
\begin{equation}  \label{eq_sumofexponent}
(\bar{A_{0}})_{i,j}=\left\{ 
\begin{array}{rl}
-a_{0}\left(\bar{x}\right) & \mbox{if $i=j$} ,\\ 
0 & \mbox{if $i \not =j$} ,
\end{array}
\right. \,, \quad
(\bar{A_{r}})_{i,j}=\left\{ 
\begin{array}{rl}
a_{r}(\bar{x}) & \mbox{if $i-j=d_{r}$}, \\ 
0 & \mbox{if $i-j \not =d_{r}$},
\end{array}
\right.
\end{equation}
where the arguments of all propensity functions are the current state $\bar{x}$ \cite{Azam_Sandu_ApproximateCME}.
Therefore the proposal distribution used in our method is:
\begin{equation}  \label{eq_proposal}
g  = \exp\left(\tau\, \sum_{r=0}^M \bar A_r\right) \cdot \delta_{\mathcal{I}(\bar{x})}\,.
\end{equation}
%

\subsection{Markov process}
\label{sect:markov}
The Markov process starts with the values of species at the current time. The candidate sample is generated by the tau-leap method. Both the candidate sample and the  current sample are evaluated based on the acceptance ratio  \eqref{eq-ratio}. The target density ratio \eqref{eq-ratio1} is 
\begin{equation} \label{alpha1}
\alpha_{1}=\frac{\pi\left(x^{*}\right)}{\pi\left(x_{t-1}\right)
}= \frac{\delta_{\mathcal{I}(x^*)}\cdot\exp\left(\tau\, \sum_{r=0}^M A_r\right)\cdot \delta_{\mathcal{I}(\bar{x})}\,}
{\delta_{\mathcal{I}(x_{t-1})}^T \cdot\exp\left(\tau\, \sum_{r=0}^M A_r\right)\cdot \delta_{\mathcal{I}(\bar{x})}\,}.
\end{equation}
For the tau-leap method $x^{*}$ is generated independent of $x_{t-1}$ and vice versa. Hence the proposal density ratio \eqref{eq-ratio2} is
%
%
%
\begin{equation} \label{alpha2}
\alpha_{2}=\frac{g\left(x_{t-1}\right)}{g\left(x^{*}\right)}= \frac{\delta_{\mathcal{I}(x_{t-1})}^T\cdot\exp\left(\tau\, \sum_{r=0}^M \bar A_r\right)\cdot \delta_{\mathcal{I}(\bar{x})}\,}
{\delta_{\mathcal{I}(x^*)}^T\cdot\exp\left(\tau\, \sum_{r=0}^M \bar A_r\right)\cdot \delta_{\mathcal{I}(\bar{x})}\,}.
\end{equation}
%
%
From \eqref{alpha1} and \eqref{alpha2} the acceptance ratio $\alpha$ is:
\begin{equation} \label{alpha}
\alpha=\frac{\delta_{\mathcal{I}(x^*)}^T\cdot\exp\left(\tau\, \sum_{r=0}^M A_r\right)\cdot \delta_{\mathcal{I}(\bar{x})}\,}
{\delta_{\mathcal{I}(x_{t-1})}^T\cdot\exp\left(\tau\, \sum_{r=0}^M A_r\right)\cdot \delta_{\mathcal{I}(\bar{x})}\,} \cdot
\frac{\delta_{\mathcal{I}(x_{t-1})}^T\cdot\exp\left(\tau\, \sum_{r=0}^M \bar A_r\right)\cdot \delta_{\mathcal{I}(\bar{x})}\,}
{\delta_{\mathcal{I}(x^*)}^T\cdot\exp\left(\tau\, \sum_{r=0}^M \bar A_r\right)\cdot  \delta_{\mathcal{I}(\bar{x})}\,}.
\end{equation}
In the acceptance/rejection test, samples which have a higher density ratio will be selected as the next state and samples which have a lower density ratio will be rejected. The Markov process samples have approximately the same density as CME (SSA) even when using a large time step in the proposal (explicit tau-leap). The only drawback of this method is the cost of performing matrix exponential. In the following section we discuss several ways to reduce this computational cost.

\section{Matrix exponential}
\label{sect:mexp}
The computation of a large matrix exponential is a problem of general interest, and a multitude of approaches suited to different situations are available.  The most straightforward, and naive, approach is a direct application of the definition of the matrix exponential
\begin{equation}
\label{eqn:matexp}
\exp(\A) = \displaystyle\sum_{k=0}^{\infty} \frac{\A^k}{k!}.
\end{equation}
While this approach is guaranteed to converge if sufficiently, possibly very many, terms are used, there are substantial numerical stability problems in the case where either the norm or the dimension of $\A$ is very large \cite{expokit, pade}.  
\subsection{Rational Approximation Methods}
Several rational approximation methods have been developed to overcome the stability and speed of convergence problems posed by the direct method.  These are based on standard function approximation methods, in the case of the Pade approximation \cite{pade}, or on the approximation of complex contour integrals, in the case of CRAM \cite{cram}.  These methods are usually paired with the ``scaling and squaring'' process of Higham \cite{Higham} to further increase the stability of the computation.
\subsubsection{Pade approximation}

The Pade approximation for $\exp(\mathbf{A})$ is computed using the $(p,q)$-degree rational function:
\[
\begin{array}{lr}
P_{pq}(\mathbf{A})=[D_{pq}(\mathbf{A})]^{-1}N_{pq}(\mathbf{A}),\\
N_{pq}(\mathbf{A})=\sum_{j=0}^{p} \frac{\left(p+q-j\right)!p! }{\left(p+q\right)!j!(p-j)!}\mathbf{A}^{j},\\
D_{pq}(\mathbf{A})=\sum_{j=0}^{q} \frac{\left(p+q-j\right)!q! }{\left(p+q\right)!j!(q-j)!}(-\mathbf{A})^{j},
\end{array}
\]
which is obtained by solving the algebraic equation:
\[
\sum_{k=0}^{\infty}\frac{\mathbf{A}^k}{k!}-\frac{N_{pq}(\mathbf{A})}{D_{pq}(\mathbf{A})}=O\left(\mathbf{A}^{p+q+1}\right)
\]
in which $P_{pq}(x)$ must match the Taylor series expansion up to order $p+q$  \cite{expokit}.
 MATLAB's \texttt{expm} function makes use of thirteenth order Pade approximation with scaling and squaring \cite{pade}. 
\subsubsection{Rational approximations of integral contours}
\label{sect:integral_contour}
In the case where the spectrum of $\A$ is confined to a region near the negative real axis of the complex plane  methods based on the rational approximation of integral contours have the potential for faster convergence than the Pade approximation. 
This approach is based on constructing parabola and hyperbola contour integrals on left complex plane and is given by \cite{cram}:
\[
r(\A)=\sum_{k=1}^{N} \frac{\alpha_k}{\A-\theta_k}
\]
Where $\theta_k$ are the quadrature points from the contour and $\alpha_k$ are the weights of the quadrature rule. \cite{cram} uses parabola approach and provides the Matlab script for both rational approximation methods and the coefficients.
\subsubsection{Chebyshev Rational Approximation Method (CRAM)}
\label{sect:Cheby}
CRAM extends the idea of rational approximation of integral contours in an attempt to obtain an optimal order of convergence. The approximation is computed as
\[
r_{k,k}(\A)=\frac{p_k(\A)}{q_k(\A)}
\]
having constraint: 
\[
\sup_{\A \in \mathbb{R}} \mid r_{k,k}(\A)-e^{\A} \mid= \inf \left\{ \sup_{\A \in \mathbb{R}} \mid r_{k,k}(\A)-e^{\A} \mid \right\}.
\]
Where $p_k$ and $q_k$ are the polynomials of order $k$. The primary difficulty in making use of the CRAM method is the procurement of suitable coefficients of the polynomials $p_k$ and $q_k$. A method for obtaining these coefficients is given in \cite{coeffs}, and they are given explicitly for $k = 14$ and $k = 16$ in \cite{cram}.

\subsection{Krylov based approximation}
\label{sect:Krylov}
For our purposes we do not seek the entire solution of $\exp(\A)$, in fact we would like only a single element of the result. Krylov based approximations get us one step closer to this ideal. Where the rational approximation methods seek to approximate the entirety of equation (\ref{eqn:matexp}), Krylov based methods seek only an approximation to
the matrix-vector product  $\exp(\A)b$.

This is done by first computing the $m$-dimensional Krylov subspace 
\begin{equation*}
\K_m = \textrm{span}\left\{b,\, \A b, \dots, \A^{m-1}b\right\}
\end{equation*}
using the Arnoldi iteration to compute the $n \times m$ orthonormal basis  matrix $\V_m$ and the $m \times m$ upper Hessenberg matrix $\Hb_m$ with $m \ll n$ such that 
\begin{equation*}
\textrm{span}(\V_m) = \K_m, \qquad
\Hb_m = \V^T \A \V.
\end{equation*}
The approximation is constructed as 
\begin{equation}
\label{eqn:kryapprox}
\exp(\A)b = \V_m \V_m^T \exp(\A) \V_m\V_m^T b = \Vert b\Vert\V_m \exp(\Hb_m) e_1 
\end{equation}
where $e_1$ is the first canonical basis vector.  The small matrix exponential exponential term in (\ref{eqn:kryapprox}) can be computed using one of the rational approximation methods with scaling and squaring extremely cheaply.
The EXPOKIT software of Sidje \cite{expokit} makes use of these techniques, with some extra consideration for Markovian cases, where the approximation of $w(t) = \exp(t \A)v$ is subject to the constraint that the resulting vector is a probability vector with components in the range of $\left[ 0, \, 1 \right]$ and the sum of these components is approximately one.

\subsection{Faster approximation techniques}
\label{sect:fast_approch}
Since we seek only a single element of the matrix exponential $(\exp(\A))_{i,j}$
we propose two techniques to speed up this computation.
\subsubsection{A single element Krylov approach} 
Using equation (\ref{eqn:kryapprox}) with $b = e_j$ leads to
\begin{equation}
\label{eqn:kryoneele}
(\exp(\A))_{i,j} = e_i^T\, \exp(\A)\, e_j = (e_i^T\, \V_m)\, (\exp(\Hb_m) \,e_1).
\end{equation}
The exponential matrix entry is computed for the cost of an $m$-dimensional Pade approximation and an $m$-dimensional dot product since $(e_i ^T\V_m)$ can be computed for ``free'' by simply reading off the $i$th row of $\V_m$, and similarly $\left(\exp(\Hb_m)e_1\right)$ is just the first column of $\exp(\Hb_m)$.  This approach avoids the construction of any additional $n$-dimensional vectors or their products.

\subsubsection{Exponentiation of a selected sub-matrix}
Computing the exponential of a large matrix is expensive. When the number of species  in a reaction system is high, the dimensions of the matrix \eqref{eq_sumofexponentexact} for target probability as well as dimensions of matrix \eqref{eq_sumofexponent} for proposal probability grow quickly. For the case of $n$ species where each has a maximum $Q$ molecules the dimension of matrix will be $\left(Q+1\right)^n \times \left(Q+1\right)^n$. 

In order to reduce costs we propose to exponentiate a sub-matrix of the full matrix. The selected rows and columns contain indices of both the current state of system at $t_n$ and candidate state at $t_n+\tau$. The motivation comes from the fact that states which are far from the current and the proposed ones do not impact significantly the acceptance/rejection test of Metropolis-Hastings algorithm and can be disregarded. Numerical experiments indicate that the error in an element $(\exp(\A))_{i,j}$ computed using a properly sized sub-matrix instead of full matrix is small. 

In order to obtain the proper size of a sub-matrix for each reaction system, we use specific information from the reaction system such as propensity functions, time step and maximum number of molecules in the system. 
Recall the general tau-leap formula \cite{Gillespie_2001}.
\[
 x\left(\bar{t}+\tau\right)=x\left(\bar{t}\right)+ \sum_{j=1}^M V_{j}\, K\left(a_j\left(x\left(\bar{t}\right)\right) \tau\right).
\]
where $K\left(a_j\left(x\left(\bar{t}\right)\right) \tau\right)$ is a random number drawn from a Poisson distribution with parameter  $a_j\left(x\left(\bar{t}\right)\right) \tau$ and $V_{j}$ is the $j$-{th} column of stoichiometry matrix. 
The expected value of the jump in the number of molecules is
\begin{equation} \label{init_ges_2}
\texttt{E}[x\left(\bar{t}+\tau\right)-x\left(\bar{t}\right)] = \sum_{j=1}^M V_{j}\, a_j\left(x\left(t_0\right)\right)\, \tau.
\end{equation}
Motivated by \eqref{init_ges_2} which indicates the weighted sum of propensities, we consider the following initial estimate of the size of the sub-matrix (S):
\begin{equation} \label{init_ges}
S \propto \frac{\|V\|}{N} \sum_{j=1}^M  a_j\left(x\left(t_0\right)\right)\, \tau \propto \bar{a} \left(x\left(t_0\right)\right)\, \tau.
\end{equation}
where $\bar{a} \left(x\left(t_0\right)\right)\ $ is the average over the propensity functions of the initial values of species. 

We seek to select a range of state indices that covers the current and proposed states. The sub-matrices are built by selecting only the rows and columns in this range from \eqref{eq_sumofexponentexact} and \eqref{eq_sumofexponent}. If the range of indices is small then the exponential computations are fast. However, if this range does not cover the representative states (both the current sample and the proposed sample), the probability ratio of the proposed sample can be far from the target probability, and the proposed sample is likely to be rejected. Choosing the size of the sub-matrix for maximum efficiency has to balance the cost of obtaining a sample (smaller is better for the cost of exponentiation) with the likelihood of accepting samples (larger is better for accuracy of approximation).
 
 

\section{Numerical experiments}
\label{sect:experiment}
This section discusses the application of the Metropolis Hastings algorithm to generate samples from the SSA distribution for three test systems:  
Schlogl \cite{Cao_2004_stability}, reversible isomer \cite{Cao_2004_stability}, and Lotka Volterra reactions \cite{Gillespie_1977}.

\subsection{Schlogl reaction}
\label{sect:Schlogl}
We first consider the Schlogl reaction system from \cite{Cao_2004_stability}

\begin{equation}
\label{eqn:schlogl}
\begin{array}{lr}
\ce{
B_{1} + 2x  <=>[\ce{c_1}][\ce{c_2}]  3x
},\\
\ce{
B_{2} <=>[\ce{c_3}][\ce{c_4}]  x
},
\end{array}
\end{equation}
whose solution has a bi-stable distribution. Let $N_1$, $N_2$ be the numbers of molecules of species $B_1$ and $B_2$, respectively.
The reaction stoichiometry matrix and the propensity functions are: 
\[
V= 
\begin{bmatrix}
1 & -1 & 1 & -1
\end{bmatrix}, \quad 
\begin{array}{l}
a_{1}(x)= \frac{c_{1}}{2}N_{1}x(x-1), \\ 
a_{2}(x) = \frac{c_{2}}{6}N_{1}x(x-1)(x-2), \\ 
a_{3}(x) = c_{3}N_{2}, \\ 
a_{4}(x) = c_{4}x. 
\end{array}
\]
The following parameter values (each in appropriate units) are used:
\[
\begin{array}{lll}
c_{1}=3 \times 10^{-7}, &c_{2}=10^{-4}, &c_{3}=10^{-3}, \\
c_{4}=3.5, &N_{1}=1 \times 10^5, &N_{2}=2 \times10^5,
\end{array}
\]
with final time $T=4$, initial conditions $x(0)=250$ molecules, and maximum values of species $Q^1=900$ molecules. 
We consider a time step $\tau=0.4$ for which the explicit tau-leap solution has a relatively large error compared to SSA. 

The initial guess for the size of sub-matrix given by \eqref{init_ges} is $250 \times 250$ and works well for the model. To accept $1,000$ samples the MCMC process
rejects about $\sim 1,200$ samples when using full matrix (whose size is $901 \times 901$). While the number of rejected using sub-matrix is approximately $1,300$. Decreasing the size of sub-matrix leads to a larger number of rejected samples. For example using a sub-matrix of size $100 \times 100$ results in approximately $2,500$ rejected samples, so this matrix size is too small. Another metric to assess whether the sub-matrix size is appropriate is the size of the residual obtained by exponentiating the full matrix and the sub-matrix. In this simulation the residual is below $10^{-8}$ for a sub-matrix size of $250 \times 250$. We have observed empirically that when the residual is larger than $10^{-2}$ the sample is likely to be rejected. The moderate number of rejected samples using the sub-matrix and the small residual indicate that the $250 \times 250$ size yields a good approximation for large matrix exponentiation.

\begin{figure}[tb]
	\begin{centering}
	\subfigure[]{
	\includegraphics[width=0.7\textwidth, height=0.45\textwidth]{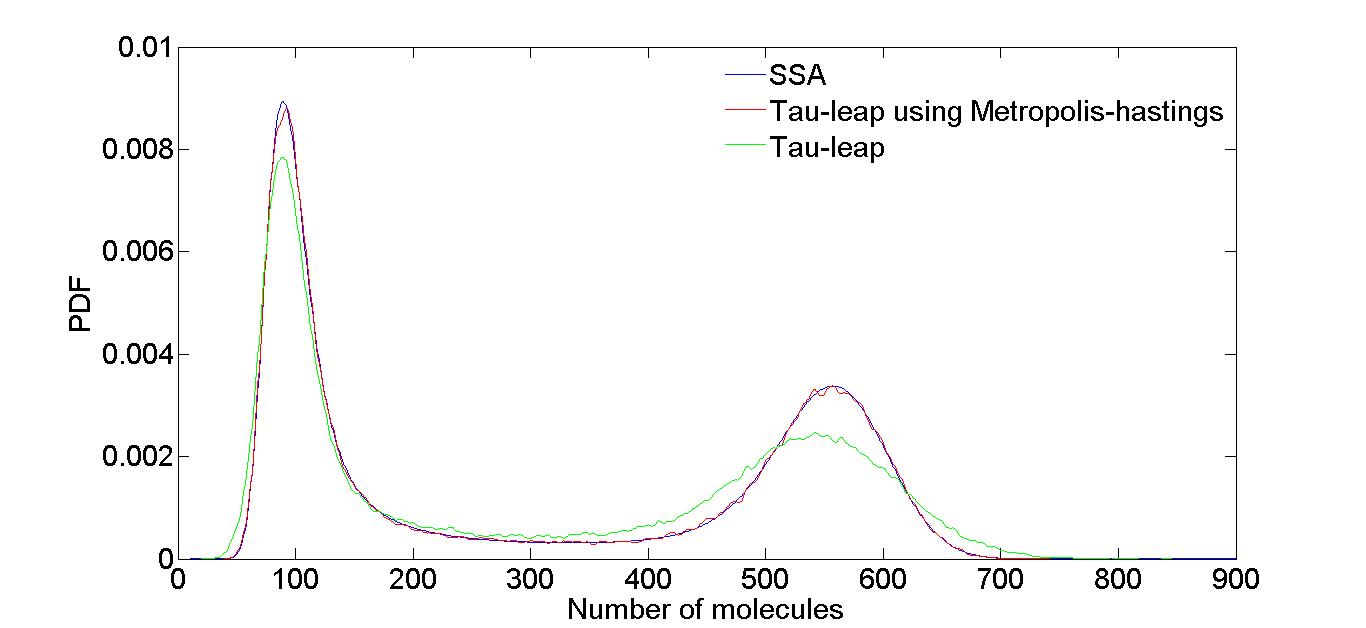} 
	\label{fig:schlogl1}
	}
	\subfigure[]{	
	\includegraphics[width=0.7\textwidth, height=0.45\textwidth]{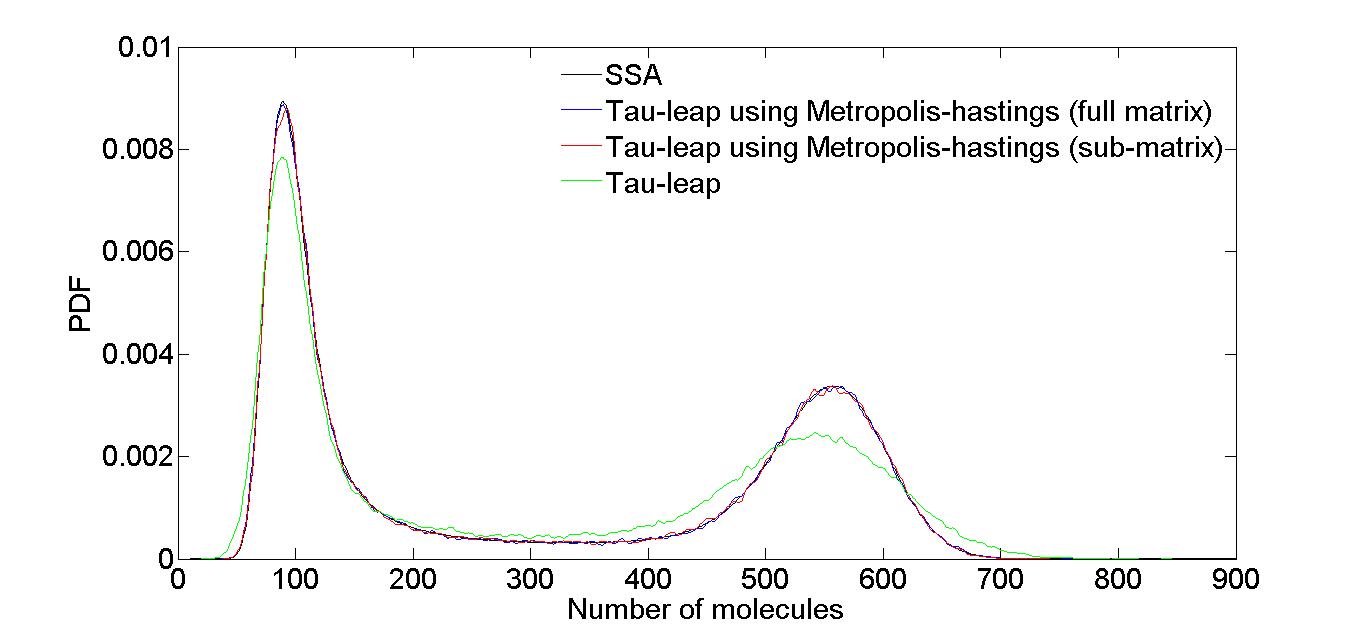} 
\label{fig:schlogl2}}
	\caption{Histograms of Schlogl system \eqref{eqn:schlogl} solutions with $\tau=0.4$ (units), final time T=4 (units), and 10,000 samples.}
	\label{fig:schlogl}
	\end{centering}
\end{figure}

%
Figure \ref{fig:schlogl1} illustrates the histogram of Schlogl reaction results obtained  by SSA, explicit tau-leap, and Metropolis Hastings using full matrix size. Figure \ref{fig:schlogl2} shows that the results obtained with a sub-matrix of size $250 \times 250$ have no visible reduction in accuracy.
Since all the eigenvalues of the matrix lie very closely to each other we employ the rational approximation technique discussed in \ref{sect:integral_contour} for exponentiating both the full matrix and the sub-matrix. The CPU time of obtaining one sample using the sub-matrix is about (0.32 sec) half the CPU time per sample when using the full matrix (0.76 sec). The approximate time of getting one sample using SSA is 0.15 sec vs. 0.02 sec. using tau-leap.

%
%
\subsection{Isomer reaction}
\label{sect:Isomer}
The reversible isomer reaction system from \cite{Cao_2004_stability} is given by:
\begin{equation}
\label{eqn:isomer}
\ce{
x_1 <=>[\ce{c_1}][\ce{c_2}]  x_2,
}
\end{equation}
and the stoichiometry matrix and the propensity functions are: 
\[
V= 
\left[\begin{array}{rr}
-1 & 1 \\
1 & -1 
\end{array}\right]\,, \qquad
\begin{array}{l}
a_{1}(x)= c_{1}x_{1} \,, ~~\\
a_{2}(x) = c_{2}x_{2} \,.
\end{array}
\]
The reaction rate values are $ c_{1}=10$, $c_{2}=10$ (units), the time interval is $[0,T]$ with $T=10$ (time units), 
the initial conditions are $x_{1}(0)=40$, $x_{2}(0)=40$ molecules, and the maximum values of the species are $Q^1=80$ and $Q^2=80$ molecules.

The estimate give by equation \eqref{init_ges} is $20$ and since this reaction system has two species the initial guess for the size of the sub-matrix is $20^2 \times 20^2$. In order to be more conservative a sub-matrix of size $500 \times 500$ is selected. In order to accept $1,000$ samples the Markov process rejects approximately $5,000$ samples when using the full matrix (of size $6,561 \times 6,561$) , and about $8,000$ samples when using the sub-matrix. Decreasing the size of sub-matrix leads to many more rejected samples. Our empirical observations show again that when the residual is larger than $10^{-2}$ the sample is likely to be rejected. We conclude that the current sub-matrix provides a good approximation for large matrix exponentiation.

\begin{figure}[tb]
	\begin{centering}
	\subfigure[]{
	\includegraphics[width=0.7\textwidth,height=0.35\textwidth]{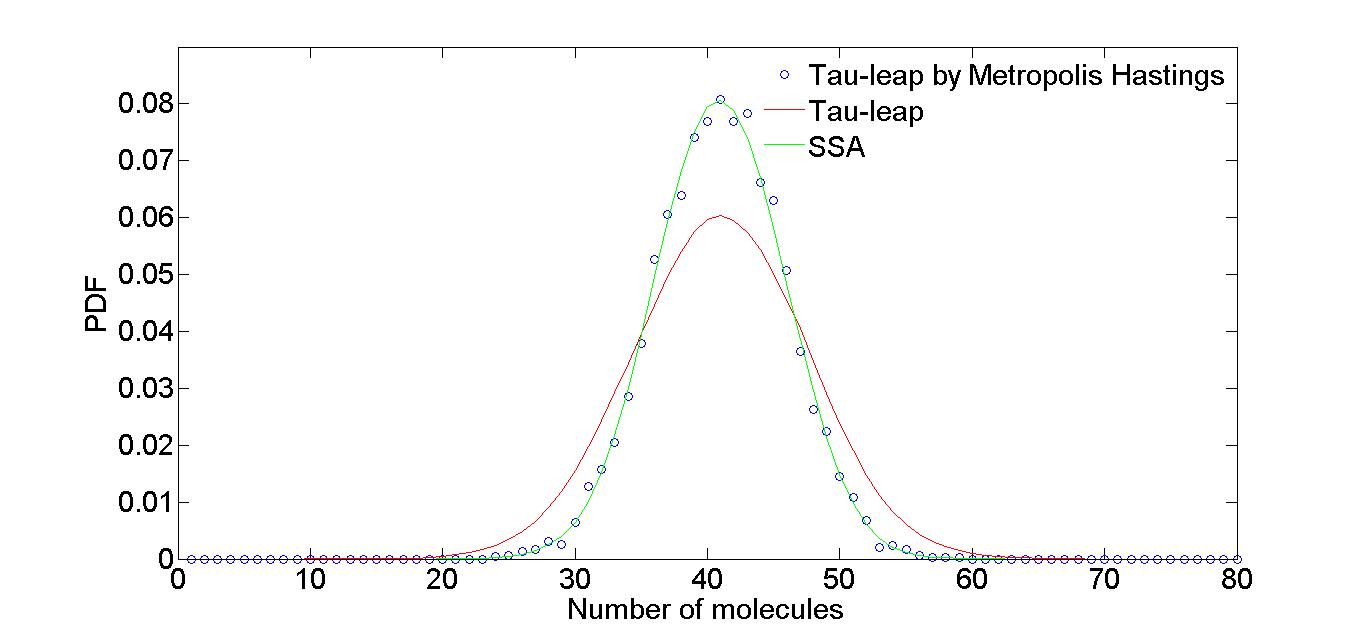} 
	\label{fig:isomer1}
	}
	\subfigure[]{	
	\includegraphics[width=0.7\textwidth,height=0.35\textwidth]{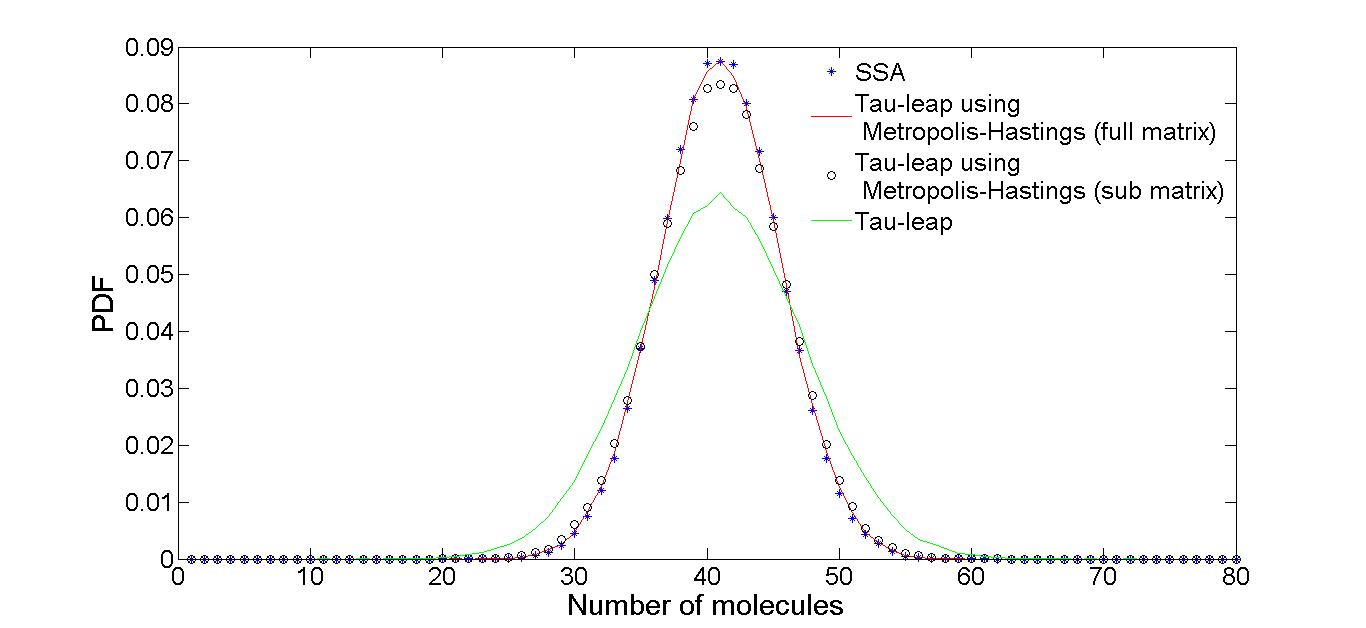} 
\label{fig:isomer2}}
	\caption{Histograms of isomer system \eqref{eqn:isomer} solutions with $\tau=0.05$ (units), final time T=1 (units), and 10,000 samples.}
	\label{fig:isomer}
	\end{centering}
\end{figure}

Figure \ref{fig:isomer1}shows the histogram of the isomer reaction solutions obtained by SSA, explicit tau-leap, and  by Metropolis Hastings using the full size matrix \eqref{eq_sumofexponentexact} and \eqref{eq_sumofexponent}. Figure \ref{fig:isomer2} shows the results using the sub-matrix of size $500 \times 500$. There is no visible reduction in accuracy. Since all the eigenvalues of matrix lie very closely to each other the exponentiation of both matrices is performed using the rational approximation technique explained in \ref{sect:integral_contour}. The CPU time of obtaining one sample using sub-matrix is 20.37 sec vs. 38.70 sec. using the full matrix. Getting a sample using SSA and tau-leap takes 0.15 sec. and 0.05 respectively.
%
%
\subsection{Lotka Volterra reaction}
\label{sect:lotka}
The last test case is Lotka Volterra reaction system \cite{Gillespie_1977}:
\begin{equation}
\label{eqn:lotka}
\begin{array}{lr}
Y+x_1 \xrightarrow {c_{1}} 2x_1,\\
x_1+x_2 \xrightarrow {c_{2}} 2x_2,\\
x_2 \xrightarrow {c_{3}} Y,\\
x_1 \xrightarrow {c_{4}} Y.
\end{array}
\end{equation}
The reaction stoichiometry matrix and the propensity functions are: 
\[
V= 
\left[ \begin{array}{rrrr}
1 & -1& 0& -1 \\[0.3em] 
0 &1 &-1&0
\end{array}\right]\,, \quad
\begin{array}{l}
a_{1}(x)= c_{1}x_{1}Y \,, ~~
a_{2}(x) = c_{2}x_{2}x_{1}, \\
a_{3}(x) = c_{3}x_{2}\,, ~~
a_{4}(x) = c_{4}x_{1}.
\end{array}
\]
The following parameter values are used (in appropriate units):
\[
\begin{array}{l}
 c_{1}=0.0002, ~~ c_{2}=0.01,~~
c_{3}=10, ~~ c_{4}=10, ~~ Y=10^{-5},
\end{array}%
\]
the final time is $T=1$, the initial conditions are $x_{1}(0)=1000$,
$x_{2}(0)=1000$ molecules, and the maximum values of species are $Q^1=2000$ and $Q^2=2000$ molecules. 
The resulting full matrix has dimension $4,004,001 \times 4,004,001$ and exponentiation is not feasible without the sub-matrix approximation.

\begin{figure}[tb]
	\begin{centering}
	\includegraphics[width=0.7\textwidth, height=0.45\textwidth]{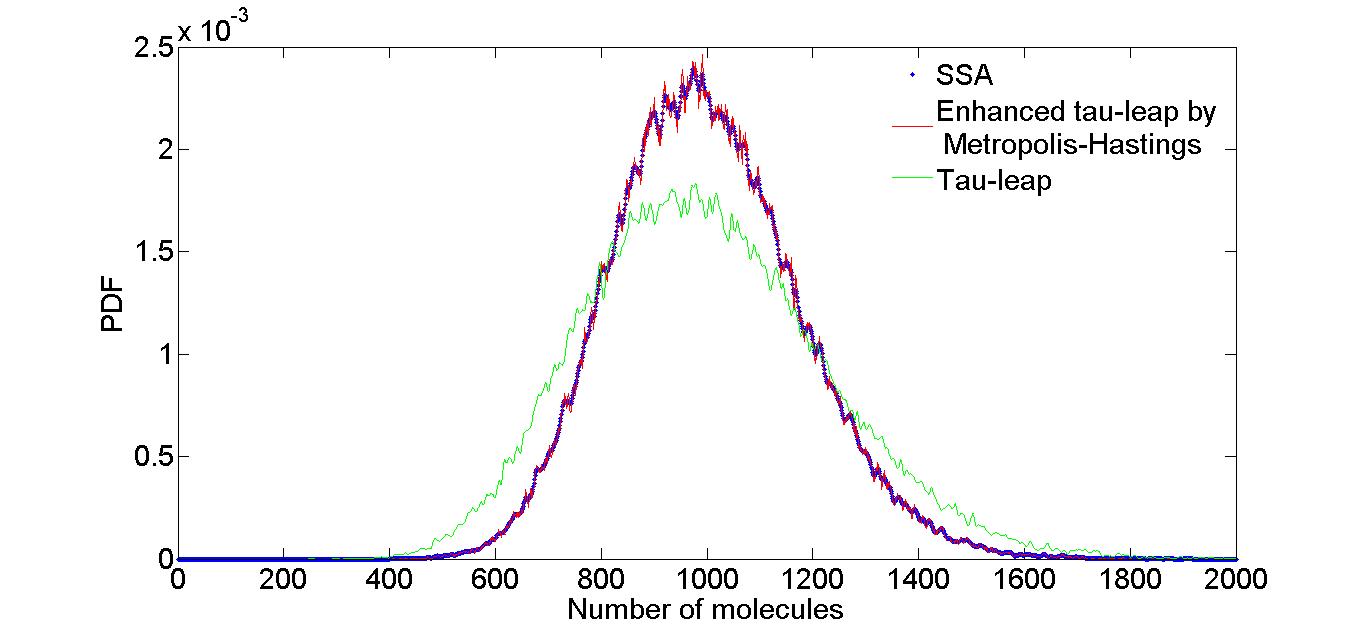} 
	\caption{Histograms of Lotka-Volterra system \eqref{eqn:lotka} solutions with $\tau=0.01$ (units), final time T=1 (units), and 10,000 samples.}
	\label{fig:lotka_res}
	\end{centering}
\end{figure}

The value predicted by equation \eqref{init_ges} is $125$, and since this reaction system has two species the initial guess for the size of sub-matrix is $15,625 \times 15,625$. This size does not work well for this system with very large number of molecules and almost never covers both current and candidate states. We increase the size of sub-matrix to  $500,000 \times 500,000$, a value obtained by trial and error. 
Figure \ref{fig:lotka_res} illustrates the histogram of Lotka-Volterra solutions obtained by SSA, tau-leap method, and Metropolis-Hastings using a sub-matrix of size discussed above. The Metropolis-Hastings sampling is very accurate. The CPU time of  matrix exponentiation using the contour integral method discussed in \ref{sect:integral_contour} is one forth of using the Krylov method stated in \ref{sect:Krylov}. However, using Krylov method gives us more accurate and stable results for large matrices than using contour integral method, hence the number of rejected samples during the Markov process using Krylov is less than the number of rejected samples using contour integral method. The CPU time of getting one sample using Metropolis-Hastings is about few hours vs. 1.51 sec. in SSA and 0.21 sec. using tau-leap. As it is clear from the numerical results one of the drawbacks of the proposed method is computing large matrix exponentiation during the Markov process. Current work of authors are focused on obtaining a faster approximation techniques to get a single element of the matrix exponential $(\exp(\A))_{i,j}$, the method which was discussed in \ref{sect:fast_approch}.

%
%

\section{Conclusions}
\label{sect:conc}
This study applies the Metropolis Hastings algorithm to stochastic simulation of chemical kinetics. The proposed approach makes use of the CME and the exponential form of its exact solution as the target probability in the Markov process. The approximation of the explicit tau-leap method is then employed for the proposal probability. The samples generated by constructing the Markov process have the same distribution as SSA even the proposals are obtained using explicit tau-leap with a large time step. Computing matrix exponentials of huge matrices can become a computational bottleneck. A practical approximation consists of selecting a sub-matrix and exponentiating it using fast approaches like Expokit and rational approximation to significantly reduce the cost of the algorithm.

\section*{Acknowledgements}
This work was partially supported by awards
NSF DMS--1419003, NSF CCF--1218454,
AFOSR FA9550--12--1--0293--DEF, AFOSR 12--2640--06,
and by the Computational Science Laboratory at Virginia Tech.

\label{sect:bib}

\bibliographystyle{plain}
\bibliography{main.bib}

\end{document}

%% file: logo.tex
\thispagestyle{empty}
\setcounter{page}{0}

\begin{Huge}
\begin{center}
Computer Science Technical Report CSTR-{\tt11} \\
\today
\end{center}
\end{Huge}
\vfil
\begin{huge}
\begin{center}
Azam S. Zavar Moosavi, Paul Tranquilli, Adrian Sandu
\end{center}
\end{huge}

\vfil
\begin{huge}
\begin{it}
\begin{center}
``{\tt Solving stochastic chemical kinetics by Metropolis-Hastings sampling}''
\end{center}
\end{it}
\end{huge}
\vfil

\begin{large}
\begin{center}
Computational Science Laboratory \\
Computer Science Department \\
Virginia Polytechnic Institute and State University \\
Blacksburg, VA 24060 \\
Phone: (540)-231-2193 \\
Fax: (540)-231-6075 \\ 
Email: \url{sandu@cs.vt.edu} \\
Web: \url{http://csl.cs.vt.edu}
\end{center}
\end{large}

\vspace*{1cm}

\begin{tabular}{ccc}
\includegraphics[width=2.5in]{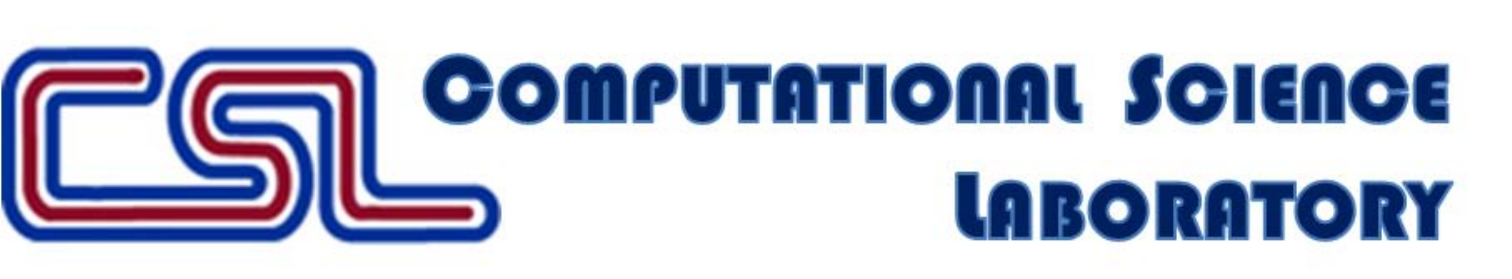}
&\hspace{2.5in}&
\includegraphics[width=2.5in]{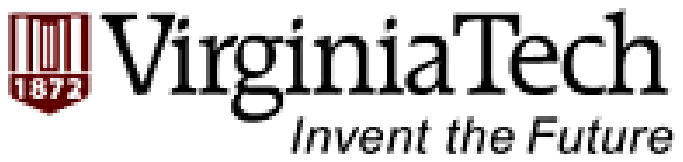} \\
{\bf\em Innovative Computational Solutions} &&\\
\end{tabular}

\newpage